\documentclass[12pt]{article}
\newcommand\f{{\bf f}}
\newcommand\g{{\bf g}}
\newcommand\h{{\bf h}}
\renewcommand\u{{\bf u}}
\renewcommand\v{{\bf v}}
\newcommand\w{{\bf w}}
\newcommand\R{{\bf R}}
\newcommand\bz{{\bf 0}}
\newcommand\hull{\mathop{\rm hull}}
\newcommand\cond{\mathop{\rm cond}}
\newenvironment{proof}{\begin{trivlist}\item[] 
{\sc Proof. }\ignorespaces}{\ \sqr\end{trivlist}} 
\newcommand\sqr{\vrule height 1.8ex width 1.0ex depth0ex} 
\newcommand\eref[1]{$(\ref{#1})$}
\newcommand\thref[1]{Theorem~$\ref{#1}$}
\newcommand\condref[1]{Condition~$\ref{#1}$}
\newcommand\sref[1]{Section~\ref{#1}}
\title{A Bernstein-B\'ezier Sufficient Condition
for Invertibility of Polynomial
Mapping Functions (Draft)}
\author{Stephen A.\ Vavasis\thanks{Supported in part by
NSF ITR Award number 0085969.}}
\newtheorem{theorem}{Theorem}
\newtheorem{condition}{Condition}
\begin{document}
\maketitle
\begin{abstract}
We propose a sufficient condition for invertibility of a polynomial
mapping function defined on a cube or simplex.  This condition is
applicable to finite element analysis using curved
meshes.
The sufficient condition is based on an analysis of the Bernstein-B\'ezier
form of the columns of the derivative.
\end{abstract}
\section{Invertibility of polynomial mapping functions}
In finite element analysis, it is common to subdivide the domain into
elements that are images of a reference domain under polynomial
functions.  This approach gives rise to the popular {\em isoparametric
elements} \cite{Johnson}.  The reference domain is 
the unit square $I^2=\{(\xi,\eta):0\le \xi,\eta\le 1\}$ or
triangle $\Delta^2=\{(\xi,\eta):0\le \xi,\eta;\xi+\eta\le 1\}$ in
$\R^2$ or the unit cube $I^3=\{(\xi,\eta,\zeta):0\le \xi,\eta,\zeta\le 1\}$
or tetrahedron $\Delta^3=\{(\xi,\eta,\zeta):0\le \xi,\eta,\zeta;\xi+\eta+
\zeta\le 1\}$ in $\R^3$.  

Polynomials defined on $\Delta^d$, $d=2,3$, generally include
monomial terms up to degree $p$ for some $p>0$.  On the other hand,
for the unit cube $I^d$, the monomial terms are generally up to degree
$p$ individually in each coordinate.  Therefore, for the rest of the
paper we use $p$ to denote the maximum total degree of polynomials on
$\Delta^2$ or $\Delta^3$, and the maximum degree in individual
coordinates of polynomials defined on $I^2$ or $I^3$.

Let $F:U\rightarrow\R^d$ be a polynomial function, where
$U$ is either $I^d$ or $\Delta^d$ and $d$ is 2 or 3.
Suppose $F$ is injective on $U$, implying that
it has an inverse $G$ defined on
$F(U)$ such that $G(F(\u))=\u$.  If such a $G$ exists and is
smooth, then we say that $F$ is {\em invertible}.
Invertibility is important for correctness of a finite element
mesh.  This definition of invertibility is sometimes called
``global invertibility'' in order to distinguish it
from
``local invertibility.''  The function $F$ is locally invertible if its
derivative, denoted as $J$, is nonsingular on the entire reference
element $U$. 
Local invertibility is necessary for global 
invertibility but is not sufficient. 
A function can be locally but not globally invertible
if $F(U)$ ``wraps around''  (e.g., consider $F(r,\theta)=r(\cos\theta,
\sin\theta)$ for $(r,\theta)\in [1,2]\times [0,4\pi]$).
Additional sufficient conditions
for global invertibility in the general setting are considered by
Ivanenko \cite{ivanenko}.
Beyond local invertibility,
one would also like an upper bound on the
``condition number'' 
\begin{equation}
\cond(F)=\max_{\u\in U}\Vert J(\u)\Vert\cdot\max_{\u\in U}\Vert J(\u)^{-1}\Vert
\label{cond}
\end{equation}
of $F$.  Even more strongly, one would
like bounds on the higher derivatives of $F$  in terms of $\cond(F)$
to satisfy the Ciarlet-Raviart \cite{Ciarlet}
conditions on convergence of order-$p$
finite element approximation.

Unfortunately, even the first problem, namely, determining nonsingularity
of the Jacobian, is a difficult problem.  We know of 
simple necessary and sufficient
conditions only for the following special cases: 
\begin{enumerate}
\item
 linear polynomials,
that is, $p=1$, defined on $\Delta^d$
of all dimensions, 
\item quadratic polynomials on $\Delta^2$
in the case $d=2$, in the special
case that $F$ is linear on two of the three edges of the reference
triangle \cite{Johnson}, and 
\item
bilinear elements (i.e., $p=1$) defined
on $I^2$ \cite{StrangFix}.
\end{enumerate}

For all other cases, we know of only separate necessary and sufficient
conditions, and this paper also proposes only a sufficient
condition.
To our knowledge,
the only necessary condition for the general case 
proposed in the literature is that
$\det(J)$ have the same sign (strictly positive or strictly negative)
at some finite list of test points.  In engineering
applications, it is common to test invertibility at the Gauss points
used for quadrature on the element \cite{Shephard}.
This test is known not to be sufficient; see further remarks on
this point below.

Sufficient conditions have been proposed for a few settings. For example
Ushakova \cite{ushakova} presents a sufficient condition for the
case $p=1$ on domain $I^3$.  Her sufficient condition is interpreted as
requiring that certain tetrahedra formed by choosing subsets of
four of the eight vertices of 
$F(I^3)$ have positive volume.  Sufficient conditions for quadratic
triangles and tetrahedra have been proposed by Salem, Canann and Saigal
\cite{SCS1, SCS2, SCS3, SCS4}.

In this work, we present a sufficient condition to ensure both
local and global invertibility
of $F$ for any of the four domains $I^d$, $\Delta^d$, $d=2,3$
and for any polynomial
degree $p$.  Our condition appears to be weaker (i.e., not able
to certify invertibility in more cases) than existing sufficient tests
for specific reference domains and values of $p$
but is considerably more general.

Our condition is based on writing $J$ in Bernstein-B\'ezier
form and then considering the convex hull of the derivatives at the
control points.  
In the next section we describe the new sufficient condition and establish
that it is sufficient for local invertibility.  In \sref{comput} we
provide an equivalent characterization of the sufficient condition that
is amenable to efficient testing.  In \sref{global} we use the equivalent
characterization to prove that the condition also implies global
invertibility.
Our condition has a desirable property that we term
``affinity,'' which we define in the last section.

\section{Bernstein-B\'ezier form}

Bernstein-B\'ezier (BB) form is a popular way to write polynomials in
computer-aided geometric design \cite{Farin}.  A univariate polynomial
of degree $p$ in BB form would be written:
\begin{equation}
F(\xi)=\sum_{i=0}^p f_i(1-\xi)^{p-i}\xi^i\frac{p!}{i!(p-i)!}
\label{BB1}
\end{equation}
and has natural parametric domain $\xi\in [0,1]$.
A bivariate polynomial with maximum degree $p$ individually
in $\xi,\eta$ is written in the form:
\begin{equation}
F(\xi,\eta)=\sum_{i=0}^p\sum_{j=0}^p \f_{i,j}(1-\xi)^{p-i}\xi^{i}
(1-\eta)^{p-j}\eta^{j}
\frac{p!p!}{i!(p-i)!j!(p-j)!}
\label{BBI2}
\end{equation}
and has as its natural parametric domain $(\xi,\eta)\in I^2$.
A trivariate polynomial with maximum degree $p$ individually
is written
\begin{equation}
F(\xi,\eta,\zeta)=\sum_{i=0}^p\sum_{j=0}^p \sum_{k=0}^p
\f_{i,j,k}(1-\xi)^{p-i}\xi^{i}
(1-\eta)^{p-j}\eta^{j}
(1-\zeta)^{p-k}\zeta^{k}
\frac{p!p!p!}{i!(p-i)!j!(p-j)!k!(p-k)!},
\label{BBI3}
\end{equation}
and has as its natural parametric domain $(\xi,\eta,\zeta)\in I^3$.

A bivariate polynomial with total degree at most $p$ is written
in the BB form
\begin{equation}
F(\xi,\eta)=\sum_{i=0}^p\sum_{j=0}^{p-i} \f_{i,j}\xi^{i}
\eta^{j}(1-\xi-\eta)^{p-i-j}
\frac{p!}{i!j!(p-i-j)!}
\label{BBD2}
\end{equation}
and has natural parametric domain $(\xi,\eta)\in \Delta^2$.

Finally, a trivariate polynomial with total degree at most $p$ is
written in the BB form
\begin{equation}
F(\xi,\eta,\zeta)=\sum_{i=0}^p\sum_{j=0}^{p-i}\sum_{k=0}^{p-i-j} \f_{i,j,k}\xi^{i}
\eta^{j}\zeta^k(1-\xi-\eta-\zeta)^{p-i-j-k}
\frac{p!}{i!j!k!(p-i-j-k)!}
\label{BBD3}
\end{equation}
and has natural parametric domain $(\xi,\eta,\zeta)\in\Delta^3$.

In all five 
cases, the  vectors $f_i$ (in 1D), $\f_{i,j}\in\R^2$ (in 2D)
or $\f_{i,j,k}\in\R^3$ (in 3D) are called the {\em control points}.
A fundamental theorem about BB form is:
\begin{theorem}
Let $U$ be the natural parametric domain of BB form for the five cases
listed above.  Then $F(U)$ is contained in the convex hull of the
control points.
\end{theorem}

Farin \cite{Farin} proves this as a consequence of the deCasteljau
algorithm for evaluating $F$, but here is a sketch of a 
more direct proof.  One observes
that on the natural parametric domain, all the factors in the summations
like $\xi$, $(1-\xi-\eta)$, etc., are nonnegative.  Furthermore, one observes
that the value of all the sums, if the control points are excluded,
is exactly 1. For example, consider \eref{BB1} without control points:
$$\sum_{i=0}^p (1-\xi)^{p-i}\xi^{i}\frac{p!}{i!(p-i)!}$$
This sum is identically 1, as seen by
considering a binomial expansion
of $(\xi+(1-\xi))^p$.  Similar argument apply to \eref{BBI2}--\eref{BBD3}.
Thus, in the natural parametric domains,
the above summations may be regarded as weighted averages of the
control points, proving the theorem.

The next feature of BB form is that once a function is in BB form, 
the derivatives can be easily put into BB form.  In all cases, the
derivative is a BB expansion of one degree lower, multiplied by $p$,
and with control points that are finite differences of the
control points for $F$ in the direction of
the variable being differentiated.
Thus, for example, in the case of \eref{BBD2},
we have
$$
\frac{\partial F}{\partial\xi}=
\sum_{i=0}^{p-1}\sum_{j=0}^{p-i-1} p(\f_{i+1,j}-\f_{i,j})\xi^{i}
\eta^{j}(1-\xi-\eta)^{p-i-j-1}
\frac{(p-1)!}{i!j!(p-i-j-1)!}
$$
and
$$
\frac{\partial F}{\partial\eta}=
\sum_{i=0}^{p-1}\sum_{j=0}^{p-i-1} p(\f_{i,j+1}-\f_{i,j})\xi^{i}
\eta^{j}(1-\xi-\eta)^{p-i-j-1}
\frac{(p-1)!}{i!j!(p-i-j-1)!}.
$$
In this case, the list of $p(p+1)/2$ vectors
of the form $p(\f_{i+1,j}-\f_{i,j})$ are control points of $\partial F/\partial
\xi$, and analogously $p(\f_{i,j+1}-\f_{i,j})$ are control points
for $\partial F/\partial \eta$.
Let us denote these two lists of control points $G_\xi$ and $G_\eta$ 
respectively.  
Similar expressions hold for the other four BB forms described above.
In 3D there is a third list $G_\zeta$.
We now state our condition:

\begin{condition}
In the case $d=2$, the matrix $[\u,\v]$ is invertible for 
every $\u\in\hull(G_\xi)$ and every $\v\in\hull(G_\eta)$.
In the case $d=3$, the matrix $[\u,\v,\w]$ is invertible
for every
$\u\in\hull(G_\xi)$, $\v\in\hull(G_\eta)$, $\w\in\hull(G_\xi)$.
\label{condit}
\end{condition}
In this condition, ``hull'' denotes the convex hull.
The main result of this section is as follows.

\begin{theorem}
If \condref{condit} holds, then the
matrix $J$ is invertible on the entire reference element.
\label{suff1}
\end{theorem}

The proof of this theorem follows from the arguments in the previous
paragraphs: \condref{condit} is sufficient since the actual Jacobians that
occur on the domain have their first columns chosen from $\hull(G_\xi)$,
etc.

\section{A computational characterization of the sufficient condition}
\label{comput}

\condref{condit} specifies our sufficient condition in somewhat
nonconstructive terms.  
In this section we provide a computational means to verify the
sufficient condition.

\begin{theorem}
In the case of $\R^2$, \condref{condit} is equivalent
to the following condition:
\begin{itemize} 
\item
there exists a vector $\h\in\R^2$ such that
for all $\f\in G_\xi\cup G_\eta$, $\h^T\f>0$, and
\item
there exists a vector $\bar\h\in\R^2$ such that
for all $\f\in G_\xi$, $\bar\h^T\f>0$ and for all $\f\in G_\eta$,
$\bar\h^T\f<0$.
\end{itemize}
\label{sep}
\end{theorem}

\begin{proof}
Let $\u,\v$ be arbitrary in $\hull(G_\xi)$ and $\hull(G_\eta)$
respectively.  This is equivalent to saying
there exist nonnegative $\alpha_1,\ldots,\alpha_n$ summing
to 1 such that $\u=\alpha_1\f_1+\cdots+\alpha_n\f_n$, where
$\f_1,\ldots,\f_n$ is an enumeration of $G_\xi$.  Similarly,
there exist nonnegative
$\beta_1,\ldots,\beta_m$ summing to 1 such that
$\v=\beta_1\g_1+\cdots+\beta_m\g_m$, where $\g_1,\ldots,\g_m$
is an enumeration of $G_\eta$.  The condition that $[\u,\v]$ is
invertible is equivalent to saying that $\u,\v$ are independent,
i.e., that there do not exist $\delta,\gamma$, at least one nonzero such that
$\delta\u+\gamma\v=\bz$.  Suppose they are dependent, and let
$\delta,\gamma$ be the two coefficients of dependence.
Without loss of generality, $\delta\ge 0$.
There are now two cases: either $\gamma\ge 0$ or $\gamma <0$.
If $\gamma\ge 0$, define $\bar\alpha_1,\cdots,\bar\alpha_n$
to be $\delta\alpha_1,\ldots,\delta\alpha_n$ and define
$\bar\beta_1,\cdots,\bar\beta_m$ to be
$\gamma\beta_1,\cdots,\gamma\beta_n$.  Then the condition that
the $\alpha$'s are nonnegative and the assumption that 
$\delta\ge 0$ is equivalent to the hypothesis that
$\bar\alpha_1,\ldots\bar\alpha_n$ are nonnegative (with no
restriction on their sum).  Similarly, the condition on the
$\bar\beta_i$'s is that they are nonnegative.

Thus, in the case that $\gamma\ge 0$, a dependence is equivalent to
the existence of nonnegative 
$\bar\alpha_i$'s and $\bar\beta_i$'s, not all zeros, such
that 
$$\sum_{i=1}^n\bar\alpha_i\f_i+\sum_{i=1}^m\bar\beta_i\g_i=\bz.$$
This problem is solved by linear programming.  By Farkas' lemma
\cite{ziegler}, there is a nonnegative solution to this problem, with
not all the coefficients zero, if and only if there does not exist a
vector $\h$ such that $\h^T\f_i>0$ and $\h^T\g_i>0$ for all
$\f_i$'s and $\g_i$'s.

Now we turn to the case of dependence when $\gamma <0$.  In this case,
using the analogous argument, there is no dependence provided that there
is a vector $\bar\h$
such that $\bar\h^T\f_i>0$ for all $i$ and $\bar\h^T(-\g_i)>0$ for all $i$.
\end{proof}

This shows that in the case of $\R^2$, 
\condref{condit} can be tested by solving two linear programming problems
over $\R^2$ to find $\h$ and $\bar\h$.  In fact, there is a much simpler
algorithm, namely, compute the args of all
the points in $G_\xi\cup G_\eta$.
To determine whether $\h$ exists, one checks whether the max arg differs
by less than $\pi$ from the min arg.  
(The branch cut for defining the arg function
 must be chosen outside the min-max range).
A similar test can determine whether $\bar\h$ exists.

In the three-dimensional setting, the algorithm is not so simple but it
is still linear time.  Following the same approach as in the previous
proof, we see that in 3D, \condref{condit} is equivalent to
the following conditions:
\begin{itemize}
\item
there exists a vector $\h_1\in\R^3$ such that
$\h_1^T\f>0$ for any $\f\in G_\xi\cup G_\eta\cup G_\zeta$,
\item
there exists a vector $\h_2\in\R^3$ such that
$\h_2^T\f>0$ for any $\f\in G_\xi\cup G_\eta\cup(-G_\zeta)$,
where $-G_\zeta$ means $\{-\f:\f\in G_\zeta\}$,
\item
there exists a vector $\h_3\in\R^3$ such that
$\h_3^T\f>0$ for any $\f\in G_\xi\cup (-G_\eta)\cup G_\zeta$, and
\item
there exists a vector $\h_4\in\R^3$ such that
$\h_4^T\f>0$ for any $\f\in G_\xi\cup (-G_\eta)\cup (-G_\zeta)$.
\end{itemize}
Each of these can be verified with linear programming, which is linear
time in three dimensions \cite{Megiddo}.  Alternatively, they can
be verified in $O(n\log n)$ time using a convex hull algorithm 
\cite{Edelsbrunner}.

\section{Global invertibility}
\label{global}

We are now in a position to prove that \condref{condit}
in fact implies global invertibility of the mapping function.
To demonstrate global invertibility, we will show that $F$ is injective, i.e.,
for all $\u,\v$ in the reference domain, 
$\u\ne\v\Rightarrow F(\u)\ne F(\v)$.  Injectivity is
sufficient for global invertibility: if an injective
smooth function has a nonsingular derivative at all points of its
domain, then the global inverse is also smooth by the inverse
mapping theorem.  The nonsingularity of the derivative is already established.

\begin{theorem}
If \condref{condit} holds, then $F$
is globally invertible on the reference element.
\end{theorem}
\begin{proof}
Let us start with the proof of injectivity in two dimensions.
Choose $\u,\v$ in the reference element such that $\u\ne\v$.
Let $\w=\v-\u$.  Then by the definition of directional derivative,
\begin{equation}
F(\v)-F(\u)=\int_0^1 J(\u+t\w)\w\, dt
\label{fvu}
\end{equation}
Assume without loss of generality
that $w_1\ge 0$, where $(w_1,w_2)$ denotes the entries of $\w$.
(The case  $w_1<0$ is handled by exchanging
the roles of $\u,\v$.) Now there are two cases: either $w_2\ge 0$
or $w_2<0$.  Suppose first that $w_2\ge 0$.  Because
we assume \condref{condit} and because we have established
\thref{sep}, we conclude that there exists a vector
$\h \in \R^2$ such that
$\h^T\f>0$ for all
$\f\in G_\xi\cup G_\eta$.  
By convexity, this implies
$\h^T\f>0$ for all
$\f\in \hull(G_\xi\cup G_\eta)$.  
This implies that $\h^TJ(\xi,\eta)$ 
is a vector both of whose
entries are positive for any $\xi,\eta$ in the reference
domain since both columns of $J$ are taken from $\hull(G_\xi\cup G_\eta)$.
Since $\w$ has both entries nonnegative, and at least one
entry of $\w$ is 
positive (because $\u\ne\v$), this means $\h^TJ(\xi,\eta)\w>0$
for all $\xi,\eta$ in the reference domain.  Thus, taking the inner
product of both sides of \eref{fvu} with $\h$ 
shows that $\h^T(F(\v)-F(\u))>0$,
and in particular, $F(\v)\ne F(\u)$.

The second case is $w_2<0$.  The argument is analogous, except we use
$\bar\h$ instead of $\h$.

Finally, the extension to 3D follows the same lines.  We assume
$w_1\ge 0$.  Depending on
the signs of entries $w_2,w_3$ of $\w$, we select one of
$\h_1,\h_2,\h_3,\h_4$.
\end{proof}

\section{Discussion}

We say that a sufficient condition for invertibility has the ``affinity''
property provided that it is always satisfied if the mapping function is
a nondegenerate affine linear mapping, or a sufficiently small perturbation
of a nondegenerate affine linear mapping.

Affinity is a desirable property since the affine linear mapping is obviously
the easiest case for invertibility.  Furthermore, if a sufficient
condition for invertibility has the affinity property,
then it can be turned into a necessary and sufficient
condition if applied to each subcell of a sufficiently fine subdivision
of the original domain.  
This is because an invertible mapping on a 
fine subdivision will behave like an affine mapping plus a small perturbation
on each subcell, and the small perturbation tends to zero as the
subdivision gets finer.

It is fairly easy to see that our test has the affinity property.  
Recall that the BB control points for a constant function are all identical.
Thus, if $F$ is affine linear, then $G_\xi$, $G_\eta$ and $G_\zeta$ are
all singleton sets.  Assuming that $F$ is affine and invertible means
that the condition must be satisfied.  Furthermore, a sufficiently small
perturbation of $F$ will give $G_\xi$, $G_\eta$, and $G_\zeta$ a positive
radius, but the condition of separability by planes will still hold.

Ushakova's condition also has the affinity
property. We suggest that any proposed sufficient condition should 
possess the affinity property in order to be considered useful.

Ushakova's sufficient condition for the $p=1$ and $U=I^3$
case is less restrictive (i.e., is able to validate more elements) than
the condition proposed here.

We conclude with a few open questions raised by this work.
\begin{enumerate}
\item
In the case $p=1$ and $U=I^3$, our condition is more restrictive
than Ushakova's.  On the other hand, it is not obvious how to generalize
her condition to higher degrees or to tetrahedral elements.  It would
be interesting to combine her techniques with the ones in this paper
to come up with less restrictive conditions for higher degrees.
Ushakova (private communication) found that experiments with
randomly generated elements indicate that for $p=1$, $U=I^3$, the conditions
stated here overly restrictive in the sense that most of the
valid (invertible)
elements among a randomly generated test set
would not satisfy the sufficient condition 
proposed here.

\item
In practical application of
isoparametric mesh generation, one defines the polynomial function only
on one face of the boundary of the element (namely, the face adjacent
to a curved exterior boundary).  The remaining coefficients defining $F$
can
be determined by a formula such as the formula proposed by Lenoir
\cite{lenoir}, which has certain theoretical guarantees.
How does the sufficient condition here specialize if
we assume that interior degrees of freedom of $F$ are determined by
Lenoir's formula?  Our very
preliminary computational tests seem to indicate that
using Lenoir's formula seems to make the element more amenable to
our sufficient condition (i.e., it seems to
minimize the diameters of the sets $G_\xi, G_\eta, G_\zeta$
compared to other choices for interior degrees of freedom).

\item
The test proposed here could be strengthened to obtain an upper bound
on $\cond(F)$ defined by \eref{cond}
in terms of the numerical values of $\h^T\f_i$ as $\f_i$ ranges
over the sets $G_\xi,G_\eta,G_\zeta$  and $\h$ is the vector
defining one of the halfspaces in \thref{sep}.
This raises the possibility of a sufficient condition not only for
invertibility but also for confirming
the Ciarlet-Raviart conditions.  Recall that the Ciarlet-Raviart conditions
require that higher derivatives of $F$ are bounded in terms of 
$\cond(F)$.

Note that other condition numbers besides $\cond(F)$ may be useful in
practice.  For example, Branets and Garanzha \cite{BG} have developed
a distortion measure related to structured grid generation.  It would
be interesting to get bounds on all of these condition numbers in terms
of the the numerical values of $\h^T\f_i$.

\end{enumerate}
\section{Acknowledgments}
The author is grateful to O.\ Ushakova, S.\ Ivanenko and V.\ Garanzha
for helpful comments on this manuscript.
\bibliographystyle{plain}
\bibliography{invert}

\begin{thebibliography}{10}

\bibitem{BG}
L.~V. Branets and V.~A. Garanzha.
\newblock Global condition number of trilinear mappings. {Application} to 3-d
  grid generation.
\newblock See {\tt http: // www.ccas.ru / gridgen}, 2001.

\bibitem{Ciarlet}
P.~Ciarlet and P.-A. Raviart.
\newblock Interpolation theory over curved elements, with applications to
  finite element methods.
\newblock {\em Computer Methods in Applied Mechanics and Engineering},
  1:217--249, 1972.

\bibitem{Edelsbrunner}
H.~Edelsbrunner.
\newblock {\em Algorithms in combinatorial geometry}.
\newblock Springer-Verlag, New York, 1987.

\bibitem{Farin}
Gerald Farin.
\newblock {\em Curves and Surfaces for Computer-Aided Geometric Design}.
\newblock Academic Press, 4th edition, 1997.

\bibitem{ivanenko}
S.~Ivanenko.
\newblock Harmonic mapping.
\newblock In J.~Thompson, B.~Soni, and N.~Weatherill, editors, {\em Handbook of
  Grid Generation}. CRC Press, 1999.

\bibitem{Johnson}
C.~Johnson.
\newblock {\em Numerical Solution of Partial Differential Equations by the
  Finite Element Method}.
\newblock Cambridge University Press, 1987.

\bibitem{lenoir}
M.~Lenoir.
\newblock Optimal isoparametric finite elements and error estimates for domains
  involving curved boundaries.
\newblock {\em SIAM J. Numer. Anal.}, 23:562--580, 1986.

\bibitem{Megiddo}
N.~Megiddo.
\newblock Linear-time algorithms for linear programming in {$R^3$} and related
  problems.
\newblock {\em SIAM J. Computing}, 12:759--776, 1983.

\bibitem{SCS1}
A.~Salem, S.~Canann, and S.~Saigal.
\newblock Robust quality metric for quadratic triangular {2D} finite elements.
\newblock In {\em Trends in unstructured mesh generation}, volume 220, pages
  73--80. AMD/ASME, 1997.

\bibitem{SCS3}
A.~Salem, S.~Canann, and S.~Saigal.
\newblock Mid-node admissible spaces for quadratic triangular {2D} arbitrarily
  curved finite elements.
\newblock {\em International J. Numer. Meth. Eng.}, 50:253--272, 2001.

\bibitem{SCS2}
A.~Salem, S.~Canann, and S.~Saigal.
\newblock Mid-node admissible spaces for quadratic triangular {2D} finite
  elements with one edge curved.
\newblock {\em International J. Numer. Meth. Eng.}, 50:181--197, 2001.

\bibitem{SCS4}
A.~Salem, S.~Saigal, and S.~Canann.
\newblock Mid-node admissible space for {3D} quadratic tetrahedral finite
  elements.
\newblock {\em Engineering with Computers}, 17:39--54, 2001.

\bibitem{Shephard}
M.~Shephard, S.~Dey, and M.~Georges.
\newblock Automatic meshing of curved three-dimensional domains: Curving finite
  elements and curvature-based mesh control.
\newblock In I.~Babuska, J.~Flaherty, J.~Hopcroft, W.\ Henshaw, J.~Oliger, and
  T.~Tezduyar, editors, {\em Modeling, Mesh Generation, and Adaptive Numerical
  Methods for Partial Differential Equations}. Springer Verlag, 1995.

\bibitem{StrangFix}
G.~Strang and G.~Fix.
\newblock {\em An Analysis of the Finite Element Method}.
\newblock Prentice Hall, 1973.

\bibitem{ushakova}
O.~Ushakova.
\newblock Nondegeneracy criteria for 3-d grid cells. {Formulas} for a cell
  volume.
\newblock SIAM J. Sci. Comp., to appear.

\bibitem{ziegler}
G.~M. Ziegler.
\newblock {\em Lectures on polytopes}.
\newblock Springer-Verlag, New York, 1995.

\end{thebibliography}
\end{document}